\documentclass[10pt,twoside]{article}
\usepackage{graphicx}
\usepackage{amsmath}
\usepackage{Latex-document}

\title{\bf  Rapid Mixing in Markov Chains\vskip 6mm}
\author{R. Kannan\vspace*{-0.5cm}\thanks{Department of Computer Science, Yale University, New Haven, CT 06520,
USA. E-mail: Kannan@cs.yale.edu}}
\date{\vspace{-8mm}}

\def\sup{\hbox{sup}}
\def\Ent{\hbox{Ent}}
\def\Min{\hbox{Min}}
\def\vol{\hbox{Vol}}

\def\prob{\hbox{ {\bf Pr}}}
\def\Min{\hbox{Min}}

\def\vol{{\rm Vol}}
\def\for{\hbox{  for  }}

\newcounter{foo}
\newcommand{\foodot}{\renewcommand{\thefoo}{\arabic{foo}.}}

\newtheorem{theorem}[foo]{Theorem\foodot}

\begin{document}

\maketitle

\thispagestyle{first} \setcounter{page}{673}

\begin{abstract}

\vskip 3mm

A wide class of ``counting'' problems have been studied in Computer Science. Three typical examples are the
estimation of - (i) the permanent of an $n\times n$ 0-1 matrix, (ii) the partition function of certain $n-$
particle Statistical Mechanics systems and (iii) the volume of an $n-$ dimensional convex set. These problems can
be reduced to sampling from the steady state distribution of implicitly defined Markov Chains with exponential (in
$n$) number of states. The focus of this talk is the proof that such Markov Chains converge to the steady state
fast (in time polynomial in $n$).

A combinatorial quantity called conductance is used for this purpose. There are other techniques as well which we
briefly outline. We then illustrate on the three examples and briefly mention other examples.

\vskip 4.5mm


\end{abstract}

\vskip 12mm

\section{Examples}\label{intro}

\vskip-5mm \hspace{5mm}

We consider ``counting problems'', where there is an implicitly
defined finite set $X$ and one wishes to compute exactly or
approximately $|X|$. In many situations, the approximate counting
problem can be reduced to the problem of generating uniformly at
random an element of $X$ (the random generation problem). This is
often the relatively easier part. Then, the generation problem is
solved by devising a Markov Chain with set of states $X$ with
uniform steady state probabilities and then showing that this
chain ``mixes rapidly'' - i.e., is close to the steady state
distribution after a number of steps which is bounded above by a
polynomial in the length of the input. [The proof of rapid mixing
is often the challenging part.] We will illustrate the problem
settings and scope of the area by means of three examples in this
section. Then we will outline some tools for proving rapid mixing
and describe very briefly how the tools are applied in some
examples. This paper presents a cross-section of methods and
results from the area. A more comprehensive survey can be found in
\cite{jerrum1}.

Our first example is the permanent of a $n\times n$ matrix $A$.
Valiant \cite{key:valiant} showed that the {\it exact} computation
of the permanent is \# P - hard, i.e., every problem in a class of
problems called \# P is reducible to the exact computation of the
permanent of a matrix; thus it is conjectured that it is not
solvable in polynomial time. The hardness result holds even for
the case with each entry a 0 or a 1 whence the problem is to find
$|X|$ where $X=\{ \sigma \in S_n : A_{i,\sigma (i)}=1 \forall
i\}.$ Note that $X$ here is implicitly defined by $A$. In the
general case, we may think of $A$ as specifying a weight $\prod_i
A_{i,\sigma (i)}$ on each $\sigma$ in $X$.

As usual, we measure running time as a function of $n$, a natural
parameter of the problem (like the $n$ above) and $1/\epsilon$,
where $\epsilon >0$ is the relative error allowed. Our primary aim
is a {\bf polynomial } (in $n,1/\epsilon$) {\bf time} bounded
algorithm; but, we will also discuss methods which help improve
the polynomial. A recent breakthrough due to Jerrum, Sinclair,
Vigoda \cite{JSV} gives an approximation algorithm with such a
time bound for the permanent (of a matrix with non-negative
entries) settling this important open problem.

Our second example starts with the classical problem of computing
the volume of a compact convex set in Euclidean $n-$ space ${\bf
R}^n$. Dyer, Frieze and the author \cite{key:dfk} gave polynomial
time algorithm for estimating the volume to any specified relative
error $\epsilon$. They first reduce the problem to that of drawing
a random point from the convex set (with uniform probability
density). They then impose a grid on space and do a ``coordinate
random walk'' - from current grid point $x$ in $K$, pick one of
the $2n$ coordinate neighbours $y$ of $x$ at random and go to $y$
if $y\in K$; otherwise, stay at $x$. Under mild conditions, it is
easy to show that the steady state distribution is uniform (over
the grid points in the set); they show that in a polynomial (in
$n$) number of steps, we are ``close'' to the steady state. [The
number of states of the chain can be exponentially large.]

Lov\'asz and Simonovits \cite{key:ls} have devised a continuous
state space random walk called the ``ball walk'' which performs
better. In this, we choose at the outset a ``step size''
 $\delta>0$. From the current point $x$, we pick at random
(with uniform density) a point $y$ in a ball of radius $\delta$
with $x$ as center. We go to $y$ if it is in $K$, otherwise, we
stay at $x$.

More generally, we may consider the integration (a ``continuous''
analog of counting) of a function over a convex set $K$. Of
particular interest are logarithmically-concave (a positive real
valued function $F$ is log-concave over a domain if $\log F$ is
concave over the domain) functions, since many families of
familiar probability density functions like the multi-variate
normal are log- concave. One may use the Metropolis version of the
random walks for convex set (cf section \ref{intro}). Rapid mixing
has been proved for this general case too \cite{AK}.

Our third set of examples concerns the Ising model and other
Statistical Mechanics problems. (see \cite{key:ising} and
references there). The computational problem arising from the
Ising model is the following : we are given a real symmetric
$n\times n$ matrix $V$ (the entries of $V$ arise as pairwise
interaction energies), a real number $B$ (the external field) and
a positive real number $\beta $ (the temperature). The Ising
partition function is defined as
$$Z=Z(V_{ij},B,\beta ) = \sum_{\sigma\in \{-1,+1\}^n } e^{-\beta H(\sigma )},
\quad\hbox{where  } H(\sigma ) = -\sum_{i,j}
V_{ij}\sigma_i\sigma_j - B\sum_k \sigma_k.$$ Jerrum and Sinclair
\cite{key:ising} presented a polynomial time approximation
algorithm to compute $Z$ in the case when all $V_{ij}$ are
non-negative (called the ferromagnetic case).

Their algorithm for the ferromagnetic case first reduces the
problem to the corresponding sampling problem and then more
interestingly reduces this sampling problem to another one where
we are given a graph (explicitly) $G(V,E)$ with positive edge
weights $w (e)$. The problem is to pick a subset of edges of $G$
at random such that the probability of picking a particular subset
$T$ is proportional to
$$w(T) = \mu^{|{\hbox{odd}(T)}|}\prod_{e\in T}w (e),$$
where $\mu$ is a given positive number and odd$(T)$ denotes the
set of odd degree vertices in $T$. [Note that in this case $X$ is
the set of all subsets of edges and we have probabilities ${\cal
P}$ on $X$ as given above, where $X, {\cal P}$ are implicitly
defined by giving $G,w$.]

\section{Preliminaries, eigenvalue connection}\label{eigen}

\vskip-5mm \hspace{5mm}

Most of what we say extends naturally to continuous state space
chains (where the set of states is (possibly uncountably)
infinite) under mild conditions of measurability, but for ease of
notation, here we state it for chains with a finite number of
states. If $P$ is the transition probability matrix with $P_{xy}$
denoting the probability of transition from state $x$ to state
$y$, for any natural number $t$, the matrix power $P^t$ denotes
the $t-$step transition probabilities. All our chains will be
connected and aperiodic and thus have steady state probabilities -
$\pi (y) = \lim_{t\rightarrow \infty }P^t_{x,y}$. ($\pi (y)$
exists and is independent of the start state $x$). [The notation
$\pi (\cdot )$ will be used throughout for steady state
probabilities.] We let the vector $p^{(t)}=p^{(0)}P^t$ denote the
probabilities at time $t$ where we start with the initial
distribution $p^{(0)}$. All our chains will be
``time-reversible'', i.e., $\pi(x)P_{xy}=\pi(y)P_{yx}$ will be
valid for all pairs $x,y$.

From Linear algebra, we get that the eigenvalues of $P$ are
$1=\lambda_1 > \lambda_2\geq \lambda_3\ldots \lambda_N\geq -1$
(where $N$ is the number of states). Standard techniques yield :
\begin{theorem}\label{theorem1}
For a finite time-reversible Markov Chain, with $\pi_0=\min_x \pi
(x)$, for any $t$,
$$\sum_x \left| p^{(t)}(x) -\pi (x)\right|
\leq {1\over\pi_0} \left[\max (|\lambda_2|,|\lambda_N| )\right]
^t.$$
\end{theorem}

Modifying a Markov Chain by making it stay at the current state
with probability 1/2 and move according to its transition function
with probability 1/2 ensures that $\lambda_N >0$ while only
increasing the (expected) running time by a factor of 2; so in the
maximum above, we need only consider $\lambda_2$. We call a chain
``lazy'' if $P_{xx}\geq {1\over 2}\forall x.$ We will use the
phrase {\bf mixing time} to denote the least positive real $\tau$
such that for any $p^{(0)}$, $\sum_x |p^{(\tau )}(x)-\pi (x)|\leq
1/4$. It is known \cite{aldous} that then for $t\geq \tau \log
(1/\epsilon )$, we have $\sum_x |p^{(t)}(x)-\pi (x)|\leq \epsilon
$.

If we have a time-reversible Markov Chain on a finite set of
states with transition probability matrix $P$ with steady state
probabilities $\pi (x)$ and $F$ is a positive real valued function
on the states, there is a simple modification of the chain with
steady state probabilities - $\pi (x) F(x) /\sum_y F(y)$, called
the the {\bf Metropolis} modification. It has transition
probabilities - $ P'_{xy}=P_{xy} \Min (1, {F(y)\over F(x)}) \for
x\not= y$. This construction is used in many instances including
as mentioned in the introduction for sampling according to
log-concave functions.

\section{Techniques for proving rapid mixing}\label{tools}

\subsection{Conductance}\vskip-5mm \hspace{5mm}

Alon and Milman \cite{key:alonmilman} and Sinclair and Jerrum
\cite{key:sj} related $\lambda_2$ to a combinatorial quantity
called ``conductance'' (in what may be looked on as a discrete
analog of Cheeger's inequality for manifolds). This has turned out
to be of great use in practice; often, first proofs of polynomial
time convergence use conductance.

For any two subsets $S,T$ of states, the {\it ergodic flow} from
$S$ to $T$ (denoted $Q(S,T)$) is defined as $Q(S,T)=\sum_{x\in S,
y\in T} \pi (x) P_{xy}$. The conductance $\Phi$ is defined by :
$$\Phi(S)= {Q(S,\bar S)\over \pi (S)}\qquad
\qquad \Phi = \min_{S:0< \pi (S)<3/4} \Phi (S).$$ $\Phi(S)$ is the
probability of escaping from $S$ to $\bar S$ conditioned on
starting in $S$ in the steady state; since $p^{(0)}$ may be this
distribution, it is intuitively clear that if the conductance of
any set is low, then the mixing time is high. More interestingly,
\cite{key:alonmilman} and \cite{key:sj} show also a converse.
\begin{theorem}\label{cond}
For a time-reversible, lazy, ergodic Markov chain with conductance
$\Phi$, we have
$$1-2\Phi \leq \lambda_2\leq 1-{1\over 2} \Phi^2.$$
\end{theorem}

While conductance has helped bound the mixing time for some
complicated chains (including the three examples mentioned in the
introduction), it is not a fine enough tool to give the correct
bounds for some simple chains. For example, consider the lazy
version of the random walk on the $2^n$ vertices on the $n-$ {\bf
cube}, where in each step, one picks at random one of the $n$
neighbours of the current vertex to go to. The mixing time is
known to be $O(n\log n)$. Conductance is  $\Theta(1/n)$ for this
example, yielding only a mixing time of $O(n^3)$ by Theorems
(\ref{cond}) and (\ref{theorem1}).

A striking contrast is the random walk on the vertices of the {\bf
cube truncated} by a half-space (i.e., the set of 0-1 vectors
satisfying a given linear inequality.) Morris and Sinclair
\cite{MS1} showed that the conductance of this walk is at least
$1/p(n)$ for a polynomial $p(\cdot )$.

We now discuss a recent improvement of conductance for chains with
a finite number of states from \cite{LK}, \cite{KLM}; similar
results hold for chains with infinite number of states. In
addition to measuring the ergodic flow from $S$ to $\bar S$, we
now also see if the flow is ``well-spread out'' in the sense that
we ``block'' a set $B\subseteq \bar S$, and then see if $Q(S,\bar
S\setminus B)$ is still high. We now define for $S$ with $0 < \pi
(S)\leq 3/4$,
$$\Psi(S) = \sup_{\alpha \in (0,\pi (S))}
\min_{B\subseteq \bar S;\;\; \pi (B)\leq\alpha } {\alpha \;
Q(S,\bar S\setminus B)\over \pi (S)^2}.$$ It is easy to show that
a set $B$ with $\pi (B) \leq {1\over 2}Q(S,\bar S)$, blocks at
most 1/2 of the flow from $S$ to $\bar S$, so we have $\Psi(S)
\geq {1\over 4} \Phi(S)^2$. Thus, an assertion that mixing time is
$O(\log (1/\pi_0) \min_S \Psi (S) )$ would be at least as strong a
result as we get from Theorems (\ref{cond}) and (\ref{theorem1}).
We prove a theorem which implies this assertion; indeed, instead
of taking $\min_S \Psi (S)$, the theorem takes an ``average'' of
this quantity over different set sizes. We say that $\psi
:[0,3/4]\rightarrow [0,1]$ is a ``blocking conductance function''
(b.c.f.) if (the second condition is technical)
\begin{eqnarray*}
\forall S, 0<\pi (S)\leq 3/4, \quad \Psi (S) \geq \psi (\pi
(S))\qquad\hbox{and}\qquad \psi (t) \leq 2\psi (t') \; \forall
0\leq t\leq t' \leq {4\over 3} t.
\end{eqnarray*}
\begin{theorem}\label{avcond}
If $\psi$ is a blocking conductance function of a lazy, ergodic,
time-reversible, finite Markov chain, with $\pi_0=\min_x \pi (x)$,
then, the mixing time is at most
$$500\int_{t=\pi _0}^{3/4}{1\over t\psi (t)}dt.$$
\end{theorem}

This has been used to improve the analysis of the ball walk for
convex sets in \cite{LK} and also some other examples in
\cite{MS2}. Also, \cite{BM} uses Theorem (\ref{avcond}) to argue
that the mixing time of the grid lattice, (in a fixed number of
dimensions) where some edges have failed according to a standard
percolation model is still within a constant of the mixing time of
the whole.

\subsection{Coupling}
\vskip-5mm \hspace{5mm}

Another important technique for proving rapid mixing is
``Coupling''\cite{aldous}. A {\it coupling} is a stochastic
process $(X_t,Y_t), t=0,1,2,\ldots $, where each of $\{
X_t,t=0,1,\ldots \}$ and $\{ Y_t, t=0,1,2,\ldots\}$ is marginally
a copy of the chain. [They may be mutually dependent.] So, we run
``two copies'' of the chain $(X_t,Y_t)$ in tandem. If $Y_0$ is
distributed according to $\pi$, the steady state distribution,
then, the distribution $p^{(t)}$ of $X_t$, satisfies
$$\sum_x |p^{(t)}(x)-\pi (x)|\leq \prob (X_t\not= Y_t).$$
To apply this, one must construct a coupling $(X_t,Y_t)$ for which
$X_t$ and $Y_t$ ``meet'' as fast as possible. This can prove
difficult. Path coupling introduced by Bubley and Dyer \cite{BD}
which we describe now simplifies the task quite a bit. In path
coupling, we have an underlying connected directed graph $G$ on
the set of states. ($G$ could just be the graph of the Markov
Chain.) $G$ defines distances between pairs of  states - namely
the length of the shortest path in $G$. We only need to define a
coupling of adjacent pairs of vertices, with the property that for
every pair of adjacent (in $G$) vertices $(u,v)$, the expected
distance between the next states of $u,v$ is at most $\beta < 1$.
They then show that
\begin{theorem}\label{pathcoupling}
If $D$ is the diameter (of $G$), then for any $t>0$, $\prob
(X_t\not= Y_t)\leq D\beta ^t$.
\end{theorem}

Propp and Wilson \cite{PW} have designed a method they call {\bf
Coupling from the Past}. This applies to chains with a partial
order on the set of states with a least state $\underline 0$ and a
greatest state $\underline 1$. They show that running two copies
of the Chain backwards - one from $\underline 0$ and one from
$\underline 1$ - with a coupling satisfying a certain monotonicity
condition until they ``meet'' gives us a good upper bound on the
number of steps needed to mix. We refer the reader to \cite{PW}
for details.
\subsection{Other methods}
\vskip-5mm \hspace{5mm}

 One way to prove a lower bound on
conductance for a chain with a finite set of states $X$ is to
construct a family of $|X|^2$ paths - one from each state to each
other using as edges the transitions of the Markov Chain, so that
no transition is ``overloaded'' by too many paths. We do not
supply here any more details of this technique referred to as the
method of ``canonical paths'' and used by Jerrum and Sinclair
\cite{JS}.

We may look upon the construction of these paths as routing a
multi-commodity flow through the network and apply techniques from
Network Flows. \cite{sinclair} pursues this. \cite{DiSt} uses
different measures of congestion to achieve improved results in
some cases and their methods are applied in \cite{FKP}.

Another important method is the use of logarithmic Sobolev
inequalities, where, we use (relative) entropy - $\Ent(p^{(t)})
=\sum_{x}p^{(t)}(x)\log {p^{(t)}(x)\over \pi (x)}$ as the measure
of distance. It is known that for ergodic Markov Chains, this
distance declines exponentially \cite{DiSa}; i.e., there is a
constant $\alpha\in (0,1)$ such that
$$\Ent (p^{(t)})\leq \alpha^t \Ent (p^{(0)}).$$
Note that $\Ent (p^{(0)})\leq \log (1/\pi_0)$. So, it suffices to choose $t=(\log\log {1\over\pi_0}+\log
(1/\epsilon ))$ $/(1-\alpha )$ to reduce the entropy to $\epsilon $; the dependence on $1/\pi_0$ is thus better.
But we need to determine $\alpha$ which is only known for simple chains. It is known that $\alpha > \lambda_2$, so
the most that this method could save over using something like Theorem (\ref{theorem1}) is the $\log (1/\pi_0)$
factor. \cite{houdre} and \cite{MS2} contain several comparisons between the log-Sobolev inequalities, eigenvalue
bounds and conductance. \cite{FK} uses the log-Sobolev inequality to prove better bounds on the Metropolis version
of the coordinate random walk for log-concave functions.

For the random walk on the cube a simple coupling argument, which,
moves both $X_t$ and $Y_t$ in the same coordinate, trying to make
them equal - shows that mixing time is $O(n\log n)$. Some
sophisticated Fourier Transform methods have been used to get much
more exact results here and the results are applicable in other
contexts too.

A traditional approach to sampling from a probability distribution
involves the so-called ``Stopping Rules'' \cite{aldous}, where one
specifies a rule for when to stop the Markov Chain and shows that
if we follow the rule, we sample (exactly) from the desired
distribution. \cite{ALW} contain results about the expected time
needed for certain stopping rules, which then serves as an upper
bound on the number of steps needed to converge.

We also mention two general techniques for deriving convergence
rates of a Markov Chain from the knowledge of convergence rates
for a simpler-to-analyze chain. The first one is called Comparison
and is developed in \cite{DS1}.  The second technique is called
Decomposition \cite{randall1}; here one decomposes the chain into
chains on subsets of states and derives a bound on the convergence
rate of the whole chain based on the rates for the ``sub-chains''
and the interconnections between them.

\section{Solution of sampling and counting problems}

\vskip-5mm \hspace{5mm}

{\bf PERMANENT}

We consider the permanent of a $n\times n$ 0-1 matrix $A$. We may
define a bipartite graph corresponding to the matrix. Each
$\sigma\in S_n$ with $A_{i,\sigma (i)}=1$ for all $i$ corresponds
to a perfect matching in the graph. Let ${\cal M}$ be the set of
perfect matchings in the graph. Unfortunately, no rapidly mixing
Markov Chain with only ${\cal M}$ as the set of states is known.
Broder \cite{broder} first defined the following Markov Chain. We
also include the set of ``near- perfect'' matchings - ${\cal M'}$
(a near-perfect matching has ${n\over 2}-1$ edges, no two
incident to the same vertex). Transitions of the Markov Chain are
as follows: In any current state, $M$, we pick an edge $e=(u,v)$
of the graph uniformly at random (all edges are equally likely)
and
\begin{itemize}
\item if $M\in {\cal M}_n $ and $e\in M$, move to $M^\prime = M-e$.
\item If  $M\in {\cal M}_{n-1}$ and $u$ and $v$ are both unmatched in $M$, then
move to $M^\prime = M+e$.
\item $M\in {\cal M}_{n-1}$, $u$ is matched to $w$ in $M$ and $v$ unmatched, then move to $M^\prime = (M+e)-(u,w)$; make a symmetric move if $v$ is matched
and $u$ unmatched.
\item In all other cases, stay at $M$.
\end{itemize}
\cite{JS} showed that if $A$ is dense (each row has at least $n/2$
1's), then the chain above mixes rapidly and in addition that
$|{\cal M'}| \leq p(n) |{\cal M}|$ for a polynomial $P(\cdot )$.
Thus, rejection sampling - accept result of a run of the chain if
the result is in ${\cal M}$ yields a polynomial time sampling
procedure.

Jerrum, Sinclair and Vigoda \cite{JSV} develop an algorithm for
the general 0-1 permanent (including the non-dense case). Here is
very brief sketch of their algorithm : An edge-weighting $w$
assigns a (positive) real weight $w(e)$ to each edge. For a
matching $M$ $w(M)=\prod_{e\in M} w(e)$ is its weight. For a set
$S$ of matchings, $w(S)=\sum_{M\in S} w(M)$. Finally, for each
pair of vertices $(u,v)$, define $w'(u,v)$ to be the ratio of the
weight of all perfect matchings to the weight of all near-perfect
matchings which leave $u,v$ unmatched. Then define the ``modified
weight'' $w'(M)$ of a matching $M$ to be $w(M)$ if $M$ is perfect
and $w(M)w'(u,v)$ if $M$ leaves $u,v$ unmatched. They first show
that a Metropolis version of the above random walk to sample
according to $w'(M)$ mixes rapidly. But the $w'$ are not known;
they argue that if we start with the complete graph and go through
a sequence of graphs, where in each step, we lower the edge weight
of a non-edge of $G$ by a factor, then we can successively
estimate $w'$ for each edge-weighting (of the complete graph) in
the sequence. The final element of the sequence has low enough
weights for the non-edges that it gives a good approximation to
the permanent.

{\bf THE ISING MODEL}

Recall the subgraph sampling problem in section \ref{intro} Here
is the random walk they use. The states of the Markov Chain are
the subsets of $E$. Their chain is the Metropolis version of the
following simple Markov Chain whose steady state probabilities are
uniform over all subsets of the edges, namely : at any current
subset $T$ of $E$, pick uniformly at random an edge $e\in E$; if
$e\in T$, then go to $T^\prime =T-e$, otherwise go to $T^\prime =
T+e$. They also make the chain lazy. The proof of a lower bound on
conductance relies on a canonical paths argument.

The algorithm that is preferred by physicists is the one due to
Swendsen and Wang \cite{key:swendsen}. This algorithm switches the
signs on large blocks of vertices of the graph at once. But while
this seems to work well in practice, no proof of rapid mixing is
known.

{\bf CONVEX SETS, LOG-CONCAVE FUNCTIONS}

Consider the ball walk in a convex set $K$ in ${\bf R}^n$ with
balls of radius $\delta $. We use the notation $P_{xy}$ for the
transition probability density from $x$ to $y$ here. The
conductance of a (measurable) subset $S$ of $K$ is now defined as
$${\int_{x\in S} \int_{y\in K\setminus S}
\pi (x) P_{xy}\over \min (\pi (S), 1-\pi(S) )}.$$ Let $\partial S$
be the boundary of $S$ interior to $K$. Since points $x\in S$ on
or near $\partial S$, intuitively have a high $\int_y P_{xy}$, a
lower bound on $\vol_{n-1}(\partial S)$ would seem to imply a
lower bound on conductance. This is indeed the case. Lower bounds
on $\vol_{n-1}(\partial S)$ have been the subject of much effort.
The most general result known is the following.
\begin{theorem}\label{isothm}
{\bf Isoperimetry} Suppose $K$ is a compact convex set in ${\bf
R}^n$ of diameter $d$ and $F$ is a positive real-valued
log-concave function on $K$. Then for any measurable $S\subseteq
K$ with $\int_S F\leq (1/2)\int_K F$, and measurable boundary
$\partial S$ interior to $K$, we have
$$\int_{\partial S} F \geq {2\over d} \int_S F.$$
\end{theorem}

The theorem was first proved for the case $F\equiv 1$ by Lov\'asz
and Simonovits \cite{key:ls} and independently also by Khachiyan
and Karzanov \cite{key:kk}. The result was generalized to the case
of general log-concave measures $F$ by Applegate and Kannan
\cite{AK} using the same techniques. We may add an extra factor of
$\ln (\int_K F /\int_S F)$ to the right hand side; this was proved
independently in \cite{LK} and also by Bobkov \cite{bobkov}. The
most recent algorithm for computing the volume of convex sets is
in \cite{KLS}, where references to earlier papers may be found.

{\bf OTHER EXAMPLES}

There are many other counting problems on which progress has been
made using this method. Again, we are not able to present a
comprehensive review here.

A notable result is the one for the truncated cube already
mentioned in section \ref{tools}. Another example of interest is
Contingency Tables - where we are given $m,n$ (positive integers)
and the row and column sums of an $m\times n$ matrix $A$. The
problem is to sample uniformly at random from the set of $m\times
n$ matrices with non-negative integer entries with these row and
column sums. The problem remains open, but there are several
partial results \cite{DKM},\cite{CD}.

There are many {\bf tiling problems}, where the problem is to pick
a random tiling of say a large square in the plane by dominoes of
a given shape. These problems arise in Statistical Mechanics. For
regular shapes, it is often possible to devise a polynomial time
algorithm to count the number exactly. But it is important to
devise algorithms with low polynomial time bounds. There has been
much progress here - see \cite{LRS} and references there. Random
generation of colorings and independent sets of a graph has
received much attention lately due to connections to Statistical
Mechanics \cite{BDGJ}.

\noindent{\bf Acknowledgment.} I thank Ravi Montenegro for
suggesting some changes in the manuscript.

\label{lastpage}

\end{document}